
\baselineskip=14pt
\parskip=10pt

\font\eightrm=cmr8 

\magnification=\magstephalf
\def\W{{\cal W}}

\def\1{{\overline{1}}}
\def\2{{\overline{2}}}
\parindent=0pt
\overfullrule=0in

\def\frac#1#2{{#1 \over #2}}

\bf
\centerline
{
The Method(!) of GUESS and CHECK
}

\rm
\bigskip
\centerline
{\it By Shalosh B. EKHAD and Doron ZEILBERGER}
\bigskip
\qquad 
{\it 
``Let us imagine the following example: A writes series of numbers down; B watches him and tries to find a {\bf law} for the sequence
of numbers. If he succeeds he exclaims: ``Now I can go on!''- So this capacity, this understanding, is something that
makes its appearance in a {\bf moment} $\dots$''}
\medskip
\qquad\qquad\qquad\qquad-Ludwig Wittgenstein, {\it Philosophical Investigations 1}, \S 151. [emphasis added]

\bigskip

{\bf Seven-Year-Old Gauss}

We all know how, {\it allegedly}, young Carl Friedrich Gauss summed-up the first $100$ natural numbers, by using
an {\it ``ingenious''} {\bf trick}. But this trick would fail miserably if he had, instead, to sum the first
$100$ perfect squares. A much better way would have been for him  to compute the {\bf first few terms}
of the sequence $a(n):=1+2+\dots+n$:
$$
a(0)=0 \quad , \quad 
a(1)=1 \quad , \quad 
a(2)=3 \quad , \quad 
a(3)=6 \quad , \quad 
a(4)=10 \quad , \quad 
$$
and then using his genius for {\it detecting patterns}, notice that
$$
a(0)=(0 \cdot 1)/2 \quad , \quad 
a(1)=( 1 \cdot 2)/2 \quad , \quad 
a(2)=( 2 \cdot 3)/2 \quad , \quad 
a(3)=(3 \cdot 4)/2  \quad , \quad 
a(4)= (4 \cdot 5)/2 \quad , \quad  
$$
then {\bf conjecture} that $a(n)=n(n+1)/2$, for every $n$. Then he could have tested it for a few more cases,
and confirmed that the conjecture is correct for $n=5$ and $n=6$, and {\it deduced}, that $a(100)=(100 \cdot 101)/2=5050$.
See [Z1]. 
[{\eightrm In fact checking it for $n=0,1,2$ constitutes a fully rigorous proof of the identity, since both sides are
polynomials of degree $2$.}]

In this case detecting a {\it pattern} was easy enough, we will soon see examples where
detecting patterns is much harder, and would need computers.

{\bf [Physical] Induction vs. Deduction}

Everyone knows that, unlike the {\it physical} sciences, that use {\it inductive logic}, and are {\it empirical},
mathematics is a {\bf deductive} science, where one proceeds from {\it axioms}, and, step by step, using
solid logic, proceeds to prove theorems. Of course, this is a gross over-simplification,
and already the great 
James Joseph Sylvester took issue with this {\it conventional wisdom} uttered by
Thomas  Huxley who claimed that

\qquad \quad ``{\it Mathematics is that study which knows nothing of experiment, nothing of induction, nothing of causation}'',

and Sylvester ([Sy]) retorted that mathematics is

\qquad \quad ``{\it  unceasingly  calling forth the faculties of observation and comparison, and one of its principal
weapons is induction, that it has frequent recourse to experimental trial and verification ...''}

And indeed, Sylvester was right that in that part of mathematics that philosophers of science call
{\it context of discovery}, mathematicians use {\it experiments} and {\it induction} (in the sense of
{\it physical} induction, not to be confused with {\it [complete] mathematical induction} that is a fully
{\it deductive} process). But in the {\it context of justification}, i.e. when presenting the {\it proof},
deductive logic rules.

{\bf Evolution Did not prepare us for reasoning logically and for rigorous proofs}

In a beautiful recent masterpiece ([Ar]), Zvi Artstein claims, very convincingly, that the reason most people
(in fact, all of them, including us, professional mathematicians!) find math so hard, is that biological evolution
did not prepare us to the rigid discipline of formal mathematical and logical reasoning. In order to
survive in the jungle, we had to use informal, intuitive, `Bayesian' `logic', if you would call it logic at all.

{\bf Mathematical Cultural Evolution Did not prepare us to Use Computers Optimally}

But in spite of that, Mathematics blossomed, and has come a long way, both as {\it queen} and {\it servant} of
the physical sciences. For more than two millenia, the {\it cultural evolution} of mathematics,
and the notions of {\it axiomatic method}  and {\it rigorous proof} ruled.
But neither Euclid, nor Gauss, not even Ramanujan, knew about the new {\it messiah}, the powerful
electronic computer, that would revolutionize  {\bf both} the {\it discovery} and the {\it justification}
of {\bf mathematical knowledge}, and would (soon!) turn
mathematics into an {\it empirical} science, just like physics, chemistry, and biology, but dealing
with mathematical entities (like numbers, equations, groups etc.) rather than with electrons and stars, (or acids and bases, or cells and genes etc.),
and we would soon abandon our fanatical insistence on {\it rigorous proof}, and
very soon {\it semi-rigorous} proofs ([Z2]) would be fully acceptable, and soon after, completely non-rigorous proofs!

We will describe below, as  {\it case studies}, many examples where `rigorous proofs' are possible, 
but {\bf not worth the trouble}, since they are so computer-time and computer-memory consuming.
In such cases, semi-rigorous, and even non-rigorous proofs, are good enough!

{\bf Revisiting an Old Chestnut: In how many ways can a gambler win n dollars and lose n dollars, and never be in debt?}

Consider a gambler who tosses a coin $2n$ times, and if it comes out Heads, wins a dollar, and if it comes out Tails,
loses a dollar. He is kicked out as soon as he is {\it in the red}, i.e. has negative capital.
In how many ways can he survive to $2n$ rounds, but at the end break even?

For example, denoting a win by $W$ and a loss by $L$, the following is the set
of such `good' {\it gambling histories}, with $n=3$:

$$
\{\quad WWWLLL \quad , \quad WWLWLL \quad , \quad  WWLLWL \quad , \quad WLWWLL \quad , \quad WLWLWL \quad \} \quad .
$$
Note, for example, that $WLLWWL$ is not allowed, since after the third round, the gambler owes a dollar.
How to count these {\it debt-less} sequences of Wins and Losses?

Let's first consider the easier problem where the gambler is allowed to be in debt.
The number of ways of doing it is ${{2n} \choose {n}}$ (pronounced ``$2n$ {\bf choose} $n$''), 
since they were $n$ wins and $n$ losses,
and Lady Luck had to {\bf choose} out of the total number of $2n$ coin-tosses, which of them are wins (and the remaining ones are losses).

Now what is the number of such histories where he was {\it never} in debt?
A classical, very famous result in enumerative combinatorics  (see for example, [Ai], p. 98) says that this number is
given by the ubiquitous {\it Catalan numbers} ([Sl])
$\frac{1}{n+1}{{2n} \choose {n}}$.
It follows that, assuming a {\it fair} coin, that if the gambler is allowed to get credit,
and tossed the coin $2n$ times, and broke even at the end, the probability that he
never had to borrow is exactly $\frac{1}{n+1}$. 

We know at least ten proofs of this result, some more elegant than others, but here  we will
present yet another, computer-assisted proof, that is very possibly the {\bf ugliest}.
Its advantage is that the same {\bf method} can prove, at least in principle, but very often
also in practice, {\bf any} problem of this type, where instead of the `atomic' stakes being taken from the set $\{-1,1\}$,
it is taken from any {\bf given} (finite) set of integers, $S$. For example,  $S=\{-1,-2,3\}$ (see below).
The same method extends to many (but not all) `walks' in two dimensions, where at each round he gets
$a$ dollars and $b$ euros where $[a,b]$ belongs to some given finite set of pairs $S$ (and where `getting' a negative
amount, means paying), and both his number of dollars and number of euros must always be non-negative.
The method extends to three and more currencies, but then the computations take way too long for today's computers.

In order to {\it illustrate} our `empirical', {\it guess-and-check} {\it method}, we will, as promised,
use it to give the $101$-th proof of the above-mentioned famous result that the number of
gambling histories with $n$ wins and $n$ losses, and never being in debt, is indeed given by the Catalan number
$\frac{1}{n+1}{{2n} \choose {n}}=(2n)!/(n!(n+1)!)$. While, like all our proofs, it was
{\it computer-generated}, it is sufficiently simple to be followed and understood
by humans with some patience.
Even the {\it discovery} of the proof, in this simple example,
could be done by hand, and we will do it.
This is for purely {\it pedagogical} reasons, to make the method understandable to humans, so that
they can understand how our silicon colleagues can handle much more complicated cases.

We will reprove the equivalent statement

{\bf Theorem}: Let $a(n)$ be the number of $n$-long sequences $w_1 \dots w_n$ with $w_i \in \{-1,1\}$, such that
$\sum_{j=1}^i w_j \geq 0$ for $1 \leq i \leq n-1$ and $\sum_{j=1}^n w_j = 0$.
The {\bf generating function}
$$
f(t):=\sum_{n=0}^{\infty} a(n) t^n \, 
$$
is given, {\bf explicitly}, by
$$
f(t) \, = \, {\frac {1 \, - \, \sqrt {1-4\,{t}^{2}}}{2\,{t}^{2}}} \quad .
$$

[Note that if $n$ is odd, $a(n)=0$, of course, and extracting the coefficient of $t^{2n}$ from the
right side gives $-\frac{1}{2} { {\frac{1}{2}} \choose {n+1}}(-4)^{n+1}$ that simplifies to the Catalan number 
$\frac{1}{n+1}{{2n} \choose {n}}$].

{\bf Proof}: Consider, more generally, {\it all} sequences $w_1 \dots w_n$ with $w_i \in \{-1,1\}$, with 
the property that the partial sums,  $\sum_{j=1}^i w_j \geq 0$, are all non-negative, but
without the restriction that the total sum,  $\sum_{j=1}^n w_j$, equals $0$.
Let the {\it weight} of such a sequence of length $n$, in $\{-1,1\}$, be $t^n x^{TotalSum}$,
in other words:
$$
Weight( w_1 w_2 \dots w_n ) :=t^n x^{w_1 + \dots + w_n} \quad .
$$

For example
$$
Weight(EmptySequence)=t^0 x^0 =1 \quad , \quad
Weight([1,1,-1,1,-1])=t^5 x^1 =t^5 x \quad , \quad
$$
$$
Weight([1,1,-1,1,-1,1,1,1,-1,-1])=t^{10}  x^2 \quad .
$$

Let $\W$ be the set of such sequences (including the empty sequence of length $0$, whose weight is $t^0x^0=1$),
and let $F(t,x)$ be the sum of all the weights, a certain 
{\it formal power series}\footnote{*}{\eightrm A formal power series is an infinite power series of the form $\sum_{n=0}^{\infty} a_n t^n$ where we do not worry
about (analytic) convergence, for example, $\sum_{n=0}^{\infty} n!^{n!} t^n$ is perfectly OK.
One can redo most of classical analysis via the algebra of formal power series, that ironically, 
in spite of the name `formal', is much more rigorous than the traditional analysis we learn in college,
since it only uses {\it finitistic} notions. },
in the variable $t$, whose
coefficients are polynomials in $x$. Note that the coefficient of $t^n x^i$ is the exact number of `good' sequences
$w_1 \dots w_n$ whose entries are taken from $\{-1,1\}$ where all the partial sums are
non-negative, and the total sum is $i$. In particular, Note that $f(t)$, our object of desire, is simply $F(t,0)$.

We claim that $F(t,x)$ satisfies the {\bf functional equation} 
$$
F(t,x)=1+ txF(t,x)+ tx^{-1} (F(t,x)-F(t,0)) \quad.
\eqno(FunctionalEquation)
$$

Indeed  consider any  $w_1 \dots w_n \in \W$. 

{\bf Case 1}: It is the empty sequence, (i.e. $n=0$), whose weight is $1$.

{\bf Case 2}: It ends with a $1$, i.e. $w_n=1$, then $w_1 \dots w_{n-1}$ is a legal sequence,
and all members of $\W$ of length $n$ that end with $1$ are obtained by taking sequences of length $n-1$ and appending $1$ to them.
Hence the sum of the weights of the members of $\W$ that end with $1$ is
$$
tx F(t,x) \quad ,
$$
since appending the $1$ increases both the $t$-count and the $x$-count by $1$.

{\bf Case 3}: It ends with a `$-1$', i.e. $w_n=-1$, then $w_1 \dots w_{n-1}$ is a legal sequence
{\bf but}, {\it in addition}, has the property that $\sum_{i=1}^{n-1} w_i >0$, so the
set of sequences of length $n$ that end with a `$-1$' are obtained by taking walks of length $n-1$,
{\bf except} those that add up to $0$, and appending `$-1$' to them.
In other words, we have to {\bf exclude} those sequence of length $n-1$ that sum-up to zero, whose
total weight is $F(t,0)$.
Hence, the sum of all the weights of the members of the set of good sequences that end with `$-1$' is
$$
tx^{-1} ( F(t,x)- F(t,0) ) \quad ,
$$
since appending the `$-1$' increases the $t$-count (the length of the walk)  by $1$, 
but decreases the $x$-count by $1$, and $F(t,x)-F(t,0)$ is the sum of the weights of walks whose total sum
is strictly positive.

This concludes the proof of the Functional Equation. Note that by iterating the mapping
$$
f(t,x) \rightarrow 1+ txf(t,x)+ tx^{-1} (f(t,x)-f(t,0)) \quad,
$$
$n$ times, starting with $f(t,x)=1$, we get an {\bf efficient}
way to generate the first $n+1$ coefficients, in $t$, of the formal power series $F(t,x)$.

We now   {\bf pull out of a hat}, the explicit expression
$$
G(t,x) \,: = \, {\frac {1-2\,xt-\sqrt {1-4\,{t}^{2}}}{2\, t \, \left( -x+t+t{x}^{2} \right) }} \quad ,
$$
and claim that $F(t,x)=G(t,x)$.
Indeed, dear readers, we are sure that, with some patience, and a bit of {\it masochism}, you would
be able to verify, {\it by hand}, the  {\it purely routine} , high-school algebra, calculation that
if one replaces $F(t,x)$ by $G(t,x)$ in Eq. $(FunctionalEquation)$, the two sides match.
In other words, please check that
$$
G(t,x)-( \quad 1+ txG(t,x)+ tx^{-1} (G(t,x)-G(t,0)) \quad ) \, \equiv \, 0 \quad.
\eqno(FunctionalEquation')
$$
If however, you are too lazy, you are welcome to `cheat' and copy-and-paste the following Maple code
to a Maple session:

{\tt   G:=(t,x) -> (1-2*x*t-sqrt(1-4*t**2))/(2*t*(-x+t+t*x**2)); \hfill\break
simplify(G(t,x)-1-t*x*G(t,x)-t/x*(G(t,x)-G(t,0))); } \quad,

and get {\tt 0}.

There is one sticky point remaining.
It is not immediately obvious that the $G(t,x)$ is a {\it formal power series} in $t$
with coefficients that are {\bf polynomials} in $x$, but
the readers are welcome to use the good old formula for solving a quadratic equation
($(-b \pm \sqrt{b^2-4ac})/2a$), to check that $G(t,x)$ is a solution of
the quadratic equation
$$
t \left( -x+t+t{x}^{2} \right) {G(t,x)}^{2} \, + \, \left( -1+2\,xt \right) G(t,x) \, +1 \,= \, 0 \quad .
\eqno(AlgebraicEquation)
$$
This can be rewritten as
$$
G(t,x) \, = \, 1+ 2xt \, G(t,x)+t \left( -x+t+t{x}^{2} \right) {G(t,x)}^{2}  \quad .
\eqno(AlgebraicEquation')
$$
that {\it manifestly} shows that $G(t,x)$ is indeed a formal power series in $t$ with coefficients that
are polynomials in $x$, and the first $n$ terms in its  Maclaurin expansion can be gotten
by  starting with $g(t,x)=1$ and applying the mapping
$$
g(t,x) \rightarrow 1+  2\,xt g(t,x)+t \left( -x+t+t{x}^{2} \right) {g(t,x)}^{2} \quad,
\eqno(AlgebraicEquation')
$$
$n$ times.

So we have proved that $F(t,x)$, the {\it weight-enumerator} (alias generating function) for our walks,
is given by the {\bf explicit} expression $G(t,x)$. 
{\bf Finally} to prove the theorem, just plug-in $x=0$  and get that
$$
f(t) \, = \, F(t,0) \, = \, G(t,0)
\,= \, {\frac {1-2\, \cdot 0 \cdot t-\sqrt {1-4\,{t}^{2}}}{2\, t \, \left( - 0 +t+t\cdot {0}^{2} \right) }}=
\, {\frac {1 \, - \, \sqrt {1-4\,{t}^{2}}}{2\,{t}^{2}}} \quad .
$$
QED!

{\bf Secrets from the Kitchen}

The above proof is the kind of proof that G.H. Hardy called, derogatorily, {\bf essentially verification}.
It pulled, {\it out of the hat}, the explicit expression $G(t,x)$ and then used {\it uniqueness} to
prove that $F(t,x)=G(t,x)$. How can me come-up with $G(t,x)$?

By guessing, of course! If we have a premonition that $F(t,x)$ satisfies, in addition to
the natural functional equation, that was derived by {\it combinatorial reasoning}, also 
an {\it algebraic} equation of the form
$$
A_0(t,x)+ A_1(t,x) F(t,x)+ A_2(t,x) {F(t,x)}^{2} \,= \, 0 \quad ,
\eqno(AlgebraicEquation)
$$
for {\it some} polynomials $A_0(t,x), A_1(t,x), A_2(t,x)$ {\bf to be determined}, one can
first {\it crank-out} the first, say, $20$ coefficients, in $t$, of $F(t,x)$ and use
{\it linear algebra} to {\it guess} the coefficients. But there is a much simpler way,
that in this, {\it toy example}, can even be done by hand.

The functional equation $(FunctionalEquation)$, can be rewritten (recall that $F(t,0)=f(t)$)
$$
F(t,x)=\frac{1-tx^{-1}f(t)}{1-tx-tx^{-1}} \quad,
$$
so let's try to {\it guess} a quadratic equation satisfied by $f(t)$. Since $f(t)$ only contains
even powers, let's consider instead $h(t):=f(\sqrt{t})$, and optimistically
look for {\it undetermined} coefficients, $c_1,c_2,c_3,c_4,c_5,c_6$, such that
$$
(c_1+c_2 t) + (c_3 +c_4 t) h(t) + (c_5 +c_6 t) h(t)^2 \equiv 0 \quad .
$$
Less wastefully, knowing the combinatorial origin of $f(t)$ (and hence $h(t)$), it is reasonable
to look for numbers $d_1,d_2,d_3$ such that
$$
h(t)=1+ d_1 t \, h(t)+ (d_2+ d_3 t) \, h(t)^2  \quad ,
\eqno(Hope)
$$
or equivalently
$$
h(t) \, - \, 1 \, - \,  d_1 t \, h(t) \, - \, (d_2+ d_3t) \, h(t)^2  \equiv 0 \quad ,
\eqno(Hope')
$$
for some {\bf numbers}, $d_1,d_2,d_3$ that are yet {\bf to be determined}. Of course, {\it a priori}, there is
no {\bf guarantee} that we would be successful, but {\bf let's try!}.

It is easy to crank out, either using $(Functional Equation)$, or even by {\it direct counting}, the
first few terms of $h(t)$ (alias $f(\sqrt{t})$)
$$
h(t)=1+t + 2 t^2 +5 t^3+14 t^4 + 42 t^5 + O(t^6) \quad .
$$
Now, plug-it-in into $(Hope')$:
$$
(1+t + 2 t^2 +5 t^3+14 t^4 + 42 t^5 + O(t^6) )-1
- d_1 t \, (1+t + 2 t^2 +5 t^3+14 t^4 + 42 t^5 + O(t^6))
$$
$$
-(d_2+d_3 t) \, (1+t + 2 t^2 +5 t^3+14 t^4 + 42 t^5 + O(t^6))^2
\equiv 0 \quad .
\eqno(Hope'')
$$
Expanding, and collecting coefficients of $t^i$ for $i=1$ to $i=5$ yields a system of five linear equations for
the three unknowns $d_1,d_2,d_3$, whose unique solution is
$d_1=0$, $d_2=0$, $d_3=1$. This leads us to {\it conjecture} that 
$h(t)$ satisfies the algebraic equation
$$
h(t)=1+ t h(t)^2  \quad ,
\eqno(Yea!)
$$
Replacing $h(t^2)$ by $f(t)$, we have {\bf guessed} that $f(t)$ is  a solution
of the quadratic equation
$$
f(t)=1+ t^2 f(t)^2 \quad .
$$
That is equivalent to the statement of the theorem, and that lead to the conjectured expression for $F(t,x)$
(i.e. the expression $G(t,x)$ above), that before we ``pulled out of the hat'', and lead to our {\it essentially verification} proof.

{\bf What is an Answer?}

In a classical paper ([W], see also [Z3]), guru Herb Wilf, defined what is a {\it satisfactory} answer to
a combinatorial enumeration problem that, typically, asks to enumerate a sequence of sets parametrized by one or more
discrete parameters. Obviously the answers $2^n$ for the number of subsets of a set of $n$ elements, and
$n!$ for the number of permutations of $n$ objects seem satisfactory. Traditionally, one wanted an
explicit `formula' like in the above two examples. But what is a formula?, it is an algorithm for getting
the answer, and some algorithms are better than others. In fact, there is always a `formula' for
enumerating {\it any} set,
$$
|A(n)|=\sum_{a \in A(n)} 1 \quad,
$$
just generate all the elements of $A(n)$, and then count them. Most combinatorial sets have {\it exponential} (or larger)
sizes, so this `formula' is (usually) useless! By using the modern computer-science {\bf dichotomy} of
{\it polynomial} vs. {\it exponential} growth, Wilf suggested that a good answer is a {\it polynomial time algorithm}.

Going back to our toy example, the functional equation $(FunctionalEquation)$ is already a good answer!
It is easy to see that it requires cubic time ($O(n^3)$) and quadratic memory ($O(n^2)$).
On the other hand, the `answers' that we found, the explicit expression for the generating function,
$f(t)=\, {\frac {1 \, - \, \sqrt {1-4\,{t}^{2}}}{2\,{t}^{2}}}$ 
is better, since it implies a quadratic in $N$ algorithm for
generating the first $N$ terms, and the explicit expression, $a(2n)=(2n)!/(n!(n+1)!)$ is yet better.

In our toy example, we found that $h(t)=f(\sqrt{t})$ satisfies the {\it algebraic} equation
$$
t \, h(t)^2 \, - \, h(t) \, + \, 1\, = \, 0 \quad,
$$
that turned out (in this simple case) to be quadratic, and hence solvable by radicals. Many combinatorial sequences
$c(n)$ have the property that their generating functions $C(t):=\sum_{n=0}^{\infty} c(n)t^n$
satisfy {\bf algebraic} equations of the form
$$
\sum_{i=0}^{d} p_i(t) C(t)^i =0 \quad ,
$$
(where $p_i(t)$ are some polynomials of $t$),
but not necessarily with $d \leq 4$, so generally they are not `solvable by radicals', but who cares?
Just giving the minimal algebraic equation satisfied by the generating function $C(t)$ of our
desired sequence is a {\bf good answer} in the sense of Wilf.

Sequences whose generating functions satisfy such algebraic equations are called {\it algebraic},
and it is well-known ([St][Z4]) that they form an {\it algebra}, and in particular it is
decidable whether $A=B$, so in order to prove that two sequences are identically equal, it
suffices to check them for finitely many cases, from $n=0$ to $n=N_0$, say, where $N_0$ can be
(usually) easily found a priori.

It is easy to see that sequences whose generating functions satisfy algebraic equations
satisfy a non-linear recurrence with {\it constant coefficients},
alas, requiring to remember all the past values. For example for the Calalan numbers (the coefficients of $h(t)$),
we have the non-linear recurrence
$$
a(n)=\sum_{i=0}^{n-1} a(i)a(n-1-i) \quad, \quad a(0)=1 \quad .
$$

But in this toy example, the sequence itself, i.e. the number of good gambling histories,
have an {\it explicit}, `closed-form' {\bf answer}, that is even better from the Wilfian point of view.
Namely, we have
$$
a(n)=\frac{(2n)!}{n!(n+1)!} \quad ,
$$
Since $a(n+1)/a(n)=\frac{4n+2}{n+2}$, an equivalent description of the sequence $a(n)$ is as the unique
solution of the linear recurrence equation
$$
(n+2)a(n+1)-(4n+2)a(n)=0 \quad , \quad a(0)=1 \quad .
$$
This is an example of yet another important {\it ansatz}, called  `$P$-recursive', or `discrete holonomic'.
A $P$-recursive sequence $\{c(n)\}$ is uniquely determined by a homogeneous linear recurrence equation with
polynomial coefficients
$$
\sum_{i=0}^{d} q_i(n) c(n+i) \, =0 \quad ,
$$
for some given {\it polynomials} $q_0(n), \dots, q_d(n)$, together with the $d$ initial values $c(0), \dots, c(d-1)$.

It is well-known ([St]) that this class, too, is an {\it algebra},  i.e. the sum and product of
$P$-finite sequences are again $P$-finite sequences, and it is always possible to decide whether $A=B$.
Of course, usually the recurrence is not {\it first-order}, like in this case, but it is a finite
and efficient {\it description} of the sequence, and once we know the recurrence, and the initial conditions,
it is possible to compute, its members in {\it linear time} and {\it constant memory}, since we can
`forget' about previous values, and only need to retain a `window' as large as the order, $d$, of the
recurrence.

Equivalently, a sequence, $c(n)$ is $P$-recursive if and only if its generating function,
$C(t)=\sum_{n=0}^{\infty} c(n)t^n$ satisfies a linear {\bf differential} equation with polynomial coefficients
$$
\left ( \, \sum_{i=0}^{L} r_i(t) \, \frac{d^i}{dt^i} \, \right ) C(t)  \, = \, 0  \quad ,
$$
with the appropriate $d$ initial conditions at $0$. Such power series are called ([St]) $D$-finite, or simply
(continuous) holonomic. The set of $D$-finite power series is also easily seen to be an algebra.

It is known ([St]), that every algebraic formal power series is always $D$-finite, and
hence its sequence of coefficients is $P$-recursive. The converse is false! For example the $P$-recursive
sequence $n!$ can be easily seen not to have an algebraic generating function.
The algebraic and holonomic ansatzes are nicely implemented in the Maple package {\tt gfun}
described in [SZ].

{\bf Back to Gambling}

Suppose that the gambler has a three-sided coin, marked with $\{-1,-2,3\}$. How many gambling histories are
there of length $n$ where (i) the gambler broke even but was never in debt (ii) he was never in debt
but not necessarily broke even at the end. Let's call the first sequence $a(n)$ and the second $b(n)$.

The second-named author discovered the following 

{\bf Theorem}: Let $b(n)$ be the number of sequences of length $n$, $w_1 \dots w_n$ with $w_i \in \{-1,-2,3\}$ such that
$\sum_{j=1}^i w_j \geq 0$ (for $1 \leq i \leq n$), and let $a(n)$ be the number of such sequences where in addition
$\sum_{j=1}^n w_j = 0$.

The sequence $a(n)$ satisfies a linear recurrence equation with (very complicated!) polynomial coefficients of order {\bf twenty} that may be viewed in

{\tt http://www.math.rutgers.edu/\~{}zeilberg/tokhniot/oW1D7}

and $b(n)$ satisfies a linear recurrence equation with (equally complicated) polynomial coefficients of order {\bf twenty one} that may be viewed in

{\tt http://www.math.rutgers.edu/\~{}zeilberg/tokhniot/oW1D8}  \quad \quad .

Let's describe how we discovered this theorem, by treating the general case.

{\bf The General Case}

Suppose you are given any (finite) set of integers $S$, and let $F(t,x)$ be the weight enumerator
according to the same weight, $Weight(w_1 \dots w_n):=t^nx^{w_1+ \dots + w_n}$, of the set of sequences in $S^n$
whose partial sums are always non-negative (equivalently, the number of one-dimensional walks starting at 0,
and never going to the negatives), then the same logic as before shows that $F(t,x)$ satisfies the
{\bf functional equation}
$$
F(t,x)=1+ t  \sum_{s \in S, s \geq 0} x^s F(t,x) +
t  \sum_{s \in S, s<0} x^s \left ( F(t,x)-\sum_{i=0}^{-s-1} x^i \cdot Coeff_{x^i} \, F(t,x) \right )  \quad .
\eqno(GeneralFunctionalEquation)
$$
For example, if $S=\{-1,-2,3\}$, the functional equation is
$$
F(t,x)= 1+ tx^3 F(t,x)+ t x^{-1} (F(t,x)-F(t,0))+ t x^{-2} (F(t,x)- F(t,0)- x \cdot Coeff_{x^1} F(t,x)) \quad .
$$

The same guess-and-check approach works for {\it any} finite set $S$. First 
we ask our beloved computer to use the functional equation to generate sufficiently many terms
of the formal power series 
$$
F(0,t) \quad , \quad Coeff_{x^1}F(t,x) \quad , \dots, \quad Coeff_{x^r}F(t,x) \quad ,
$$
where $r=-min(S)-1$.
Then ask the computer to {\it guess}
algebraic (or holonomic) descriptions for these, and call these conjectured representations
$g_0(t),g_1(t), \dots,  g_r(t)$. Now define $G(t,x)$ by replacing $F(0,t)$ by $g_0(t)$,
$Coeff_{x^1}F(t,x)$ by $g_1(t)$ etc., and derive an algebraic (or holonomic) representation for
$G(t,x)$. Now go backwards! Once you have an algebraic (or holonomic) description of $G(t.x)$, you also
have algebraic (or holonomic) descriptions of $G(t,0)$, $Coeff{x^1}G(t,x)$, etc., and checking that
$(GeneralFunctionalEquation)$ with $F(t,x)$ replaced by $G(t,x)$ is a {\bf purely routine}
computation in the algebraic (or holonomic) ansatz!

Since we know that we {\bf can} do it, {\bf why bother?}, let's just stick to the empirically-guessed
algebraic (or holonomic) descriptions for $f(t)$, and declare it a theorem! Also, since $G(t,x)$
satisfies {\it some} (often, very complicated!) algebraic (and differential) equation with
polynomial coefficients in $t$ and the parameter (`catalytic' variable $x$), we know that
$G(t,1)$, the generating function for all good sequences (not necessarily totaling $0$) is
algebraic (and $D$-finite, and hence its sequence of coefficients, satisfies {\it some} linear
recurrence with polynomial coefficients), it must be the same one that we guessed!
So let's {\bf declare} it a theorem, and leave the proof to those (obtuse!) readers who care about
{\it mathematical certainty}. 

It turns out that, for this one-dimensional case, there is a `meta-theorem', due to Phillippe Duchon ([D], see also [AZ])
that {\it guarantees} that everything in site is algebraic (and hence holonomic), so we know {\it a priori},
that there is {\it some} linear recurrence with polynomial coefficients satisfied by our
enumerating sequence, and if we did not find it, it is just because we didn't look far enough.
Also once found, we are {\it guaranteed} that there {\bf exists} a fully rigorous proof, but carrying out
all the details would take much longer than the initial discovery of the `conjectured' recurrences, hence
it is a waste of time to actually do the full proof.

{\bf Two Dimensions and Beyond}

For walks in two dimensions, staying in the first (non-negative) quadrant,
with a prescribed set of steps, $S$, it is no longer guaranteed that
the counting sequences are holonomic, as shown, in a seminal paper, by
Marni Mishna and Andrew Rechnitzer [MR].
Nevertheless, there are many cases where it {\bf is} holonomic, and
some very smart people (See [BM] and its many references) are trying to understand why, and whenever possible,
use human ingenuity (possibly assisted by computers), using such sophisticated
methods as the {\it kernel method} (perfected in [BM])
to {\bf derive} such equations.

Here too it is easy to set-up a functional equation in each case (but now we have {\it two}
`catalytic' variables, $x$ and $y$, in addition to $t$ (time), and the weight is
$t^n x^i y^j$ where $n$, as before is the length of the walk, and $[i,j]$ is the end-point of the walk
(that must be in the positive quadrant). If the set of steps is $S$, then the same reasoning as before
shows that $F(t,x,y)$ satisfies the functional equation
$$
F(t,x,y)= 1+ 
$$
$$
t \sum_{[s_1,s_2] \in S} 
x^{s_1} y^{s_2} \left ( F(t,x)- \sum_{i=0}^{-s_1-1} Coeff_{x^i} F(t,x,y) - \sum_{j=0}^{-s_2-1} Coeff_{y^j} F(t,x,y)
+\sum_{i=0}^{-s_1-1} \sum_{j=0}^{-s_2-1} Coeff_{x^i y^j} F(t,x,y)
\right )
$$
[Note that if $s_1 \geq 0$ or $s_2 \geq 0$ , then the corresponding sums are empty.]

In {\bf each} of the known cases, our {\bf naive} guess-and-check approach works just as well!
And with much less human effort! Not only that, we can do some cases where
the kernel method fails. We now no longer have an {\it a priori} guarantee that
the proof, in the holonomic ansatz, would succeed, but it is extremely unlikely, that
once we succeeded in finding conjectured
linear differential equations with polynomial coefficients
satisfied by $F(t,0,0)$ and by $F(t,1,1)$ that $F(t,x,y)$ would not be a solution of
another, more complicated, such linear differential equation with coefficients
that are now polynomials in $x$, $y$, and $t$. If such a differential equation,
for $F(t,x,y)$, exists, then its proof is {\bf purely routine}, but finding it 
(usually) takes too long, so why bother?

If you really want to play it safe, then try and find
$D$-finite descriptions for, say $F(t,2,0)$ and $F(t,0,2)$ (since here we only have the variable $t$, this
requires far less computational resources), and you may even venture, say $F(t,2,3)$.
If all these special cases lead to $D$-finite descriptions, then you can be  {\bf really sure}.

\vfill\eject

{\bf Maple packages and Sample Input and Output Files}

This article is accompanied by four Maple packages: {\tt W1D}, for one-dimensional walks, {\tt W1Dp}, for
the probability analog, where the die may be loaded, {\tt W2D}, for two-dimensional walks, and
{\tt W3D}, for three dimensions. They all may be gotten from the front of this article

{\tt http://www.math.rutgers.edu/\~{}zeilberg/mamarim/mamarimhtml/gac.html} \qquad ,

that contains links to many {\bf articles}, generated by these packages, about walks for many sets of steps.
For some simple cases (small sets in one dimensions), we present fully rigorous proofs, or at
least sketches, but for the more complicated cases (for example for the set $S=\{-2,-1,3\}$,
quoted above, we did not bother, since we {\bf don't care}.

{\bf Looking for Patterns}

Francis Bacon, the pioneer of the {\it Scientific method}, believed that
the scientist should study nature without any prejudice, find patterns, and then
generalize, and formulate theories. But since the haystack is so big, people
came to realize that there is always some pre-conceived ideas, and all
observation is {\it theory-laden}, and  one looks for specific types of patterns.

In the present case studies, we knew what {\it kind} of patterns to look for, 
namely  the algebraic or holonomic {\it ansatzes},
but the {\it specific} patterns were way too complicated for the
naked brain. So it usually takes much longer than a {\bf moment} (see Wittgenstein's quote in the motto), even for a powerful computer,
to detect  non-trivial patterns. It is our duty, as human coaches, to discover new ansatzes that would fit data
not accounted by already known
ansatzes, and then ask computers to search for patterns within them. One of the simplest ansatzes, that could have
been used by 7-year-old Gauss, is the {\it polynomial ansatz}, and it is used often in IQ tests, where one is
asked to continue a sequence, and a few such examples are given by Wittgenstein in his {\it Philosophical Investigations}.\footnote{**}
{\eightrm
Wittgenstein objected that one can continue a sequence arbitrarily and then make-up some `law' that justifies it.
But with the combination of Occam's razor, and having specific ansatzes in mind, and only using, say, half of the
known sequence for the `guessing', and being able to confirm our guess with the remaining known values,
resolves Wittgenstein's objection.}
But humans can only detect trivial and superficial patterns, so instead of looking for patterns themselves,
they should teach the computer to look for them.

{\bf Conclusion: It's time to Make Rigorous proofs Optional}

With all due respect to counting walks, for us, it is but a {\it case study} to illustrate
a class of problems where  fully rigorous proofs can be safely abandoned. We believe that this
would be the case, in at most fifty years, for the rest of mathematics. We will come to realize
that fully rigorous proofs are only possible for relatively trivial statements, for example
Fermat's Last Theorem, the Poincar\'e conjecture, and the Four Color Theorem.
But for really deep (and interesting!)
mathematical knowledge, we would have to be content, if lucky, with {\it semi-rigorous} proofs, where
we know that a proof exists, but it is too complicated for us, and even for our computers, to find it, and
more often with fully non-rigorous (heuristic and empirical) proofs.

{\bf References}

[Ai] Martin Aigner, {\it ``A Course in Enumeration''}, Springer, 2007.

[AZ] Arvind Ayyer and Doron Zeilberger, {\it Two dimensional directed lattice walks with boundaries},
in: {\it ``Tapas in Experimental Mathematics''} (Tewodros Amdeberhan and Victor Moll, eds.), 
Contemporary Mathematics {\bf 457} (2008), 1-20. \hfill\break
{\tt http://www.math.rutgers.edu/\~{}zeilberg/mamarim/mamarimhtml/twoDwalks.html} \qquad .

[Ar] Zvi Artstein, {\it ``Mathematics and the Real World: The Remarkable Role of Evolution in the Making of Mathematics''}, Prometheus Books, 2014. 

[BM] Mireille Bousquet-M\'elou and  Marni Mishna, {\it Walks with small steps in the quarter plane},
in: ``Algorithmic Probability and Combinatorics'', Contemporary Mathematics {\bf 520} (2010), 1-40. 
\hfill\break
http://arxiv.org/abs/0810.4387 \qquad .

[D] Phillippe Duchon, {\it On the enumeration and generation of generalized Dyck words},
Discrete Mathematics {\bf 225} (2000), 121-135. \hfill\break
{\tt www.labri.fr/perso/duchon/Papiers/Gen-Dyck.ps } \qquad .

[MR] Marni Mishna and Andrew Rechnitzer, {\it Two non-holonomic lattice walks in the quarter plane},
Theor. Computer Science {\bf 410} (38-40)(2009), 3616-3630. \hfill\break
{\tt http://arxiv.org/abs/math/0701800} \qquad .

[Sl] Neil Sloane, {\it ``The On-Line Encyclopedia of Integer Sequences''}, Sequence {\bf A000108}, \hfill\break
{\tt https://oeis.org/A000108} \qquad .

[St] Richard Stanley, {\it Differentiably finite power series}, European J. Combinatorics {\bf 1} (1980), 175-188. \hfill\break
{\tt http://www-math.mit.edu/\~{}rstan/pubs/pubfiles/45.pdf} \qquad .

[Sy] James Joseph Sylvester, {\it Inaugural Presidential Address to the Mathematical and Physical Section of the British
Association at Exter}, Aug. 1869. Reprinted in: ``The Law of Verse'', London, Green and Co., 1870, 101-130.
Also {\it Collected Works} {\bf v. 2} \#100, 650-661.

[SZ] Bruno Salvy and Paul Zimmermann, {\it GFUN: a Maple package
for the manipulation of generating and holonomic functions
in one variable}, ACM Trans. Math. Software {\bf 20} (1994), 163--177.

[W] Herbert S. Wilf, {\it What is an answer?}, American Mathematical Monthly {\bf 89}(1982), 289-292.

[Z1] Doron Zeilberger, {\it Opinion 129: The ``Lost'' Diary of Carl Friedrich Gauss Should Be Made Public}, April 1, 2013,
{\tt http://www.math.rutgers.edu/\~{}zeilberg/Opinion129.html} \qquad .

\vfill\eject

[Z2] Doron Zeilberger,  {\it Theorems for a price: tomorrow's semi-rigorous mathematical culture}, \hfill\break
Notices of the Amer. Math. Soc. {\bf 40} \# 8 (Oct. 1993), 978-981. Also Math. Intell. {\bf 16} \#4 (Fall 1994), 11-14. \hfill\break
{\tt http://www.math.rutgers.edu/\~{}zeilberg/mamarim/mamarimhtml/priced.html}   \qquad .

[Z3] Doron Zeilberger, {\it Enumerative and Algebraic Combinatorics}, in:
{\it ``Princeton Companion to Mathematics''} (W. Timothy Gowers, ed.), Princeton University Press, 2008, pp. 550-561. \hfill\break
{\tt http://www.math.rutgers.edu/\~{}zeilberg/mamarim/mamarimPDF/enu.pdf } \qquad .

[Z4] Doron Zeilberger, {\it An Enquiry Concerning Human (and Computer!) [Mathematical] Understanding},
in: {\it ``Randomness \& Complexity, from Leibniz to Chaitin''}  (C.S. Calude, ed.), World Scientific, Singapore, 2007, 383-410. \hfill\break
{\tt http://www.math.rutgers.edu/\~{}zeilberg/mamarim/mamarimhtml/enquiry.html} \qquad .

\bigskip
\bigskip
\hrule
\bigskip
Doron Zeilberger, Department of Mathematics, Rutgers University (New Brunswick), Hill Center-Busch Campus, 110 Frelinghuysen
Rd., Piscataway, NJ 08854-8019, USA. \hfill \break
zeilberg at math dot rutgers dot edu \quad ;  \quad {\tt http://www.math.rutgers.edu/\~{}zeilberg/} \quad .
\bigskip
\hrule
\bigskip
Shalosh B. Ekhad, c/o D. Zeilberger, Department of Mathematics, Rutgers University (New Brunswick), Hill Center-Busch Campus, 110 Frelinghuysen
Rd., Piscataway, NJ 08854-8019, USA.
\bigskip
\hrule

\bigskip
Published in The Personal Journal of Shalosh B. Ekhad and Doron Zeilberger  \hfill \break
({ \tt http://www.math.rutgers.edu/\~{}zeilberg/pj.html})
and {\tt arxiv.org}. [Also Submitted, by invitation of Martin Aigner, to a new edition of {\it Alles Mathematik} edited by him.]
\bigskip
\hrule
\bigskip
{\bf Feb. 15, 2015}

\end